# A learning model predictive control for virtual coupling in railroads

**Miguel A. Vaquero-Serrano**[1]  |  **Francesco Borrelli**[2]  |  **Jesus Felez**[1]

[1]Department of Mechanical Engineering, Universidad Politécnica de Madrid, Madrid, Spain

[2]Department of Mechanical Engineering, University of California at Berkeley, Berkeley, CA 94701, USA

**Correspondence**

Jesus Felez, Department of Mechanical Engineering, Universidad Politécnica de Madrid, 28016, Madrid, Spain
Email: jesus.felez@upm.es

**Funding information**

Spanish Science and Innovation Ministry – State Research Agency. Grant: PID2021-124761OB-I00

**ABSTRACT**

The objective of this paper is to present a novel intelligent train control system for virtual coupling in railroads based on a Learning Model Predictive Control (LMPC). Virtual coupling is an emerging railroad technology that reduces the distance between trains to increase the capacity of the line, whereas LMPC is an optimization-based controller that incorporates artificial intelligence methods to improve its control policies. By incorporating data from past experiences into the optimization problem, LMPC can learn unmodeled dynamics and enhance system performance while satisfying constraints. The LMPC developed in this paper is simulated and compared, in terms of energy consumption, with a general MPC, without learning capabilities. The simulations are divided into two main practical applications: a LMPC applied only to the rear trains (followers) and a LMPC applied to both the followers and the first front train of the convoy (leader). Within each application, the LMPC is independently tested for three railroad categories: metro, regional, and high-speed. The results show that the LMPC reduces energy consumption in all simulation cases while approximately maintaining speed and travel time. The effect is more pronounced in rail applications with frequent speed variations, such as metro systems, compared to high-speed rail. Future research will investigate the impact of using real-world data in place of simulated data.

## 1 | INTRODUCTION

Nowadays, railroads must confront a paramount challenge related to the congestion of the infrastructure (Sharma, Kandpal, & Santibanez Gonzalez, 2024). Instead of building new tracks, a great effort has been given to the optimization of the usage of the existing networks (Dolinayova, Zitricky, & Cerna, 2020). In this context, new opportunities arise in the form of the use of data for the implementation of intelligent solutions, by means of the artificial intelligence (AI), and in the form of new rail signaling and control systems, by means of the virtual coupling (VC). In the introduction of this article, both technologies are outlined to conclude with the main objective of this article: the proposal of a novel train controller for virtual coupling based on intelligent algorithms able to optimize the behavior of each train within a convoy.

As for the AI approaches in railroads is concerned, they can be analyzed from three different perspectives: AI techniques, AI research, and AI applications (Besinovic, et al., 2022). These three perspectives, especially those related to the techniques and applications, can be used to: i) analyze relevant recent literature in the sector, ii) understand where the proposed controller fits inside the paradigm of intelligent algorithms, and iii) compare the proposed controller with the most extended alternative solutions.

Within the AI techniques, relevant methods are summarized, such as the evolutionary computing and machine learning (ML) (Tang, et al., 2022). First, evolutionary computing comprises procedures such as genetic programming (Jeschke, Sun, Jamshidnejad, & De Schutter, 2023) and particle swarm optimization (Wang, Zhang, Wang, & Zhang, 2019). Second, ML groups several types of algorithms, recently outstanding the use of reinforcement learning (RL) algorithms. Inside ML, deep learning (DL) algorithms are generally characterized by the use of artificial neural networks (ANN), such as convolutional neural networks (CNN), recurrent neural networks (RNN) and long-short term memory (LSTM) networks (Wang, et al., 2025).

Within the AI research fields, some disciplines like expert systems, data mining, and adversarial search are introduced (Tang, et al., 2022).

Last, within the AI applications, the AI has relevant applications in computer vision and image processing, operations research, and autonomous systems (Besinovic, et al., 2022) in rail domains such as maintenance and inspection (M&I), traffic planning and management (TP&Mgt), and autonomous driving and train control (AD&TC) (Tang, et al., 2022). According to the authors, the M&I rail domain has been the traditional research topic in railroads, followed in research effort by the TP&Mgt rail topic. Despite attracting more attention in recent years, the AD&TC rail topic remains as one of the least researched fields, which might be connected to a lack of widespread publicly available data for training most of the models compared to other rail domains (Pappaterra, Flammini, Vittorini, & Bešinović, 2021). To illustrate this, some examples can be introduced regarding each rail domain.

As for the M&I rail domain is concerned, recent research has focused on the use of ANNs to recognize defects in high-speed slab tracks using limited datasets (Cai, et al., 2024), as well as on the improvement of data variety by synthesizing crack images with a generative adversarial network (GAN) in order to ultimately include them in the learning process of detection algorithms based on DL, as in (Shim, 2024) and (Huang, et al., 2024). Moreover, some works have focused on the detection of corrosion anomalies in steel bridges by means of a siamese CNN (Ghiasi, et al., 2025) and on the monitoring of vibrations by using ensemble learning (Zhuang, Liu, & Tang, 2024). Furthermore, LSTM networks have been used to analyze lateral accelerations of the car body and alert for non-programmed maintenance inspections to prevent comfort losses (Garrido Martínez-Llop, Sanz Bobi, & Olmedo Ortega, 2023).

As for the TP&Mgt rail topic is concerned, the typical addressed issues are routing, timetabling (Liu, Dabiri, Wang, & De Schutter, 2024), shunting (Ying, Chow, & Chin, 2020), and trajectory (or speed profile) generation (Li, Or, & Chan, 2023), as well as traffic analysis to supervise delays, predict conflicts or reschedule as a response to disruptions (Šemrov, Marsetič, Žura, Todorovski, & Srdic, 2016). Some examples of AI techniques applied to address these problems are random forests (Kecman & Goverde, 2015), swarm intelligence (Wang, Zhang, Wang, & Zhang, 2019), genetic algorithms (Zhou, Lu, & Wang, 2022), approximate dynamic programming (Wang, Trivella, Goverde, & Corman, 2020), and, more often, RL, such as in (Tang, Chai, Wu, Yin, & D'Ariano, 2025), (Wang, et al., 2022), (Su, et al., 2022), (Lin, et al., 2023), and (Liu, Lin, & Liu, 2024).

Final, as for the AD&TC rail topic is concerned, a subclassification may be established based on (Tang, et al., 2022), (Hewing, Wabersich, Menner, & Zeilinger, 2020), and the works of (Kim, Tay, Guanetti, & Borrelli, 2019) for road vehicles: model dynamics for the controller, model dynamics for the simulation environment, state estimator, and intelligent control. In addition, and in relation to the latter, the AI application in control can involve either substituting the controller or just improving the controller with the aim of: (*i*) enhancing its execution, or (*ii*) modeling uncertainties in order to improve the control policy through a terminal set.

Thus, regarding recent research in railroads, some works have focused on the development of dynamical models for simulation environments based on RNNs that allow the testing of classical controllers (Liu, Yang, & Yang, 2022). In addition, He, Lv, Liu, and Tang (2022) and Liu, Chai, Liu, Wang, and Chai (2022) developed some LSTM-based models to predict future states of the in-front train.

However, intelligent control has been the main research area during the last years within the AD&TC rail topic. It has focused on reducing energy consumption and improving driving efficiency, mainly through the use of RL, as in other rail topics. For example, Basile, Lui, Petrillo, and Santini (2024) use deep reinforcement learning (DRL) to replace a classical controller and to control the coordination and maneuvering of heterogeneous high-speed trains under nonlinearity and uncertainty. Nevertheless, there is an emerging AI technique that is starting to be applied in railroads: the learning model predictive control (LMPC).

LMPC is an evolution of model predictive control (MPC) that integrates learning algorithms from past experiences within the optimization problems that characterize the classical MPCs. LMPC has been successfully used with different conceptualizations in fields such as road vehicles, power systems, and smart infrastructures (Arroyo, Manna, Spiessens, & Helsen, 2022). These conceptualizations vary from including ML techniques within the optimization problem to the usage of data so as to improve the control policy through a terminal set. For instance, Arroyo, Manna, Spiessens, and Helsen (2022) includes RL in MPC for a smart energy management in buildings, whereas Kim, Tay, Guanetti, and Borrelli (2019) apply LPMC in road vehicles focusing on the terminal set approach.

Considering the above literature and according to Bertsekas (2024) and Dobriborsci, Osinenko, and Aumer (2022), LMPC has the advantage of reaching a compromise between the use of data, which characterizes the ML techniques, and the use of known parameters, such as the system dynamics and constraints, which are typical of the control theory. This compromise allows, not only the use of a reduced amount of data compared with other ML techniques, but the improvement of the previously known dynamics through learning, and always ensuring the constraints' satisfaction and, therefore, its safety and stability. Thus, due to the fact that the real-time optimization control problem has always a predominant role, LMPC is suited for applications where the satisfaction of constraints is critical, including the online control of transport vehicles. As previously mentioned, the use of mixed approaches is also possible, for example, the use of DRL in MPC, as in (Arroyo, Manna, Spiessens, & Helsen, 2022) and (Airaldi, Schutter, & Dabiri, 2023), or the use of LSTM networks in MPC (Su, Chai, Chen, & Lv, 2021). However, Bertsekas (2019) concludes that a policy may be more accurate by using specific algorithms than by using ANN-based RL.

As for railroad applications are concerned, LMPC has been recently applied to the TP&Mgt rail topic (Liu, da Silva, Dabiri, Wang, & De Schutter, 2025). However, as far as we are concerned, LMPC without ML techniques has not been applied in the AD&TC rail topic for train control applications that, unlike previous references, could be considered as an operational layer.

As for the new rail signaling and control system is concerned, virtual coupling (VC) is an emerging technology that allows an increase of the capacity of the line in comparison with other existing signaling systems under the concepts of fixed block or moving block. This increase is obtained by means of the relative braking distance concept, which generally considers, not only the speed of the controlled train, but also the speed of the in-front train in order to calculate the braking curve and, therefore, to calculate the minimum allowed distance between trains (Vaquero-Serrano & Felez, 2023a).

Due to the fact that the controlled train employs speed information about neighboring trains, VC usually distinguishes between a 'train' and a 'convoy'. Thus, a convoy is a group of $N$ trains that are not mechanically coupled, but are said instead to be virtually coupled because each one of them moves according to the relative braking distance. Moreover, within the convoy with trains named after $n = \{1, \dots, N\}$, different categories of trains may also be observed: the 'leader', which occupies the first position ($n = 1$) and is the first train within the convoy in the rolling direction, and the 'followers', which refer to the rest of the trains. This differentiation is necessary because the leader should move according to conventional signaling systems, whereas the followers will be governed by the relative braking distance, and therefore, they are the ones that run under VC conditions. In addition, given a general train $n$, the train $n-1$ is generally referred to as the 'in-front train' or 'preceding train', whereas the train $n+1$ is generally referred to as the 'rear train' or 'following train'. Nevertheless, additional considerations could me made when considering information from more than two trains, which depends on the communication topology selected in each case (Vaquero-Serrano & Felez, 2023b).

Regarding recent developments in the field of VC, Aoun, et al. (2021) analyzed the market potential of VC and Aoun, Goverde, Nardone, Quaglietta, and Vittorini (2024) compared VC and moving block signaling systems based on a risk analysis. In addition, Quaglietta, Wang, and Goverde (2020) developed a multi-state model that allows the integration of VC with conventional signaling systems. Later, this model was enhanced with the inclusion of a dynamic safety margin to consider risk factors in railway operations (Quaglietta, Spartalis, Wang, Goverde, & van Koningsbruggen, 2022).

Regarding the control of the train under VC conditions, Felez, Kim, and Borrelli (2019) presented a nominal MPC, whereas Vaquero-Serrano and Felez (2023a) presented a min-max MPC to deal with errors and uncertainties. There are also recent examples of controllers aimed at dealing with errors and uncertainties by means of a tube-MPC (Liu, Zhou, Su, Xun, & Tang, 2023) and an artificial potential field (Ji, Quaglietta, Goverde, & Ou, 2025), besides additional controllers based on cooperative control (Liu, Dabiri, Wang, Xun, & De Schutter, 2024), RL (Liu, Lang, Luo, Tang, & Chai, 2024), and digital twin-driven control based on RL (Ye, et al., 2025). Moreover, being of interest the problem of VC at low speeds and the stops at the stations (Lang, Liu, Luo, & Lin, 2022), Luo, Tang, Chai, and Liu (2024) recently developed a controller for stopping at the stations. For more information on previous VC controllers, see reviews of (Xun, Li, Liu, Li, & Liu, 2022), (Wu, Ge, Han, & Liu, 2023), and (Felez & Vaquero-Serrano, 2023).

Furthermore, there have been some early practical implementations of up to 60 km/h (Mujica, Henche, & Portilla, 2021) and up to 80 km/h (Liu, Luo, Tang, Zhang, & Chai, 2024).

Table 1 summarizes and compares the references related to intelligent control applied to VC, especially focusing on the AD&TC rail domain.

**TABLE 1** Recent references involving intelligent control for VC in some railway domains.

| Reference | AI technique | AI rail application |
|---|---|---|
| (Basile, Lui, Petrillo, & Santini, 2024) | RL | AD&TC |
| (Liu, Lang, Luo, Tang, & Chai, 2024) | RL | AD&TC |
| (Ye, et al., 2025) | RL | AD&TC |
| (Liu, da Silva, Dabiri, Wang, & De Schutter, 2025) | LMPC | TP&Mgt |
| This study | LMPC | AD&TC |

In conclusion, this article presents a novel LMPC for the emerging technology of VC in railroads. The novelty of this work comes from the fact that the LMPC that improves the control policy through a terminal set has not been applied to train control in an operational layer until now, as far as we are concerned. Therefore, the objective of the paper is to optimize the behavior of a virtually coupled train convoy, maintaining as a control policy the maximum possible speed allowed by the limitations and characteristics of the line, and keeping the components of the convoy as close as possible in safety conditions, that is, without risk of collision. Under these conditions, the effect of the LMPC is analyzed by evaluating the energy consumption as a result of the controller's learning.

The remainder of this paper is organized as follows. Section 2 describes the dynamic model that is used as the control model. Section 3 describes the formulation of the proposed LMPC. Section 4 defines the simulation cases and discusses their corresponding results. Section 5 includes the conclusions and main findings of the paper.

## 2 | DYNAMIC MODEL

In this section, we describe the dynamic model used as the control model. The presentation begins with the formulation of the dynamic equations, which are subsequently vectorized to enable a compact representation suitable for the LMPC formulation introduced in the next section. Key variables relevant to the control problem are also defined.

The dynamic equations on which the controller is based are written in Equation (1). The description of each variable, including the definition of units, can be seen in Table 2. This table also includes all the variables used in this paper.

$$\dot{s} = v \tag{1a}$$

$$\dot{v} = (F - R_T)/M \tag{1b}$$

$$\dot{F} = (u - F)/\tau \tag{1c}$$

**TABLE 2** Description of the variables used in the formulation of the controller.

| Variable | Units | Description | Variable | Units | Description |
|---|---|---|---|---|---|
| $s$ | m | Position of the head of the train | $v_{min}$ | m/s | Minimum allowable speed |
| $v$ | m/s | Speed of the train | $v_{max}$ | m/s | Maximum allowable speed in the line |
| $F$ | N | Current force applied by the actuators of the train | $j_{max}$ | m/s³ | Maximum allowable traction and braking jerk |
| $\tau$ | s | Average time constant of the actuators | $d_{des}$ | m | Desired distance between trains at a standstill |
| $M$ | kg | Train mass | $d_{min}$ | m | Minimum allowable distance between trains |
| $u$ | N | Input force, calculated as the decision variable by the controller | $u_{max}$ | N | Maximum allowable traction and braking force |
| $A$ | N | First resistance coefficient of the train resistance | $P_{max}$ | W | Maximum allowable traction and braking power |
| $B$ | N/(m/s) | Second resistance coefficient of the train resistance | $v_{DP}$ | m/s | Maximum speed at each point $s$ according to the speed profile curve calculated using a dynamic programming approach. It is used as a reference speed |
| $C$ | N/(m/s)² | Third resistance coefficient of the train resistance | $v^{\lambda}$ | m/s | Running speeds from previous iterations of the LMPC algorithm introduced as input data to the controller according to the current position $s$ |
| $g$ | m/s² | Gravitational acceleration | $Q_f$ | - | Learning costs of previous iterations for each associated $v^{\lambda}$. Introduced as input data |
| $i$ | m/m | Grade of the track. Positive values indicate an uphill grade; negative values indicate a downhill grade | $\lambda$ | - | Learning variable that weights speeds for the terminal speed. It is calculated as an internal variable by the controller |
| $R$ | m | Absolute value of the curve radius of the track | $\varepsilon^{v_{max}}$ | m/s | Slack variable for the maximum speed constraint |
| $L$ | m | Train length | $\varepsilon^{v_{min}}$ | m/s | Slack variable for the minimum speed constraint |
| $a_l$ | m/s² | Maximum braking of the preceding train | $\varepsilon^{d_{rel}}$ | m | Slack variable for the minimum relative braking distance constraint |
| $a_f$ | m/s² | Maximum braking of the train | $\varepsilon^{d}$ | m | Slack variable for the minimum distance constraint |
| $d$ | m | Absolute distance between trains | $t_s$ | s | Integration step |
| $d^{ret}$ | m | Relative braking distance between trains | $H_p$ | - | Prediction horizon |
| $j$ | m/s³ | Longitudinal jerk of the train | $H_c$ | - | Control horizon |
| $s^p$ | m | Predicted position for the preceding train | $\mathbf{X}$ | Multiple | State vector $\mathbf{X} = [s, v, F]^T$ |
| $v^p$ | m/s | Predicted speed for the preceding train | $R_T$ | N | Total train's resistance |

where:

$$R_T = A + Bv + Cv^2 + Mgi + M \cdot {}^{6}/_{R} \qquad (2)$$

Defining the state vector $\mathbf{X} = [s, v, F]^T$, the dynamic equations can be vectorized and incorporated into the state equation written in Equation (3). In this Equation, and in reference to the subscript notation $(k + 1|t)$, the state vector $\mathbf{X}$ at the time step $k + 1$ is calculated starting from the information available at time $t$, where $k$ is the integration step ranging from $t$ to the prediction horizon $(H_p)$.

$$\mathbf{X}_{k+1|t} = \mathbf{X}_{k|t} + t_s \dot{\mathbf{X}}_{k|t} \qquad (3)$$

There are also some key variables required for the LMPC formulation—namely the longitudinal jerk $(j_{k|t})$, the absolute distance between trains $(d_{k|t})$, and the relative braking distance between trains $(d_{k|t}^{rel})$, which are defined in Equations (4)-(6), respectively. These variables at the prediction time step $k$ are calculated starting from the information available at time $t$, where $k$ is the integration step ranging from $t$ to $H_p - 1$ for the jerk variable and to $H_p + 1$ for the distance variables.

$$j_{k|t} = \frac{u_{k+1|t} - u_{k|t}}{M \cdot t_s} \qquad (4)$$

$$d_{k|t} = s_{k|t}^p - s_{k|t} - L \qquad (5)$$

$$d_{k|t}^{rel} = s_{k|t}^p - s_{k|t} - L + \frac{\left(v_{k|t}^p\right)^2}{2 \cdot a_l} - \frac{\left(v_{k|t}\right)^2}{2 \cdot a_f} \qquad (6)$$

## 3 | LMPC FORMULATION

The objective of this Section is to describe the formulation of the proposed LMPC. For this purpose, this Section is divided into three subsections. Subsection 3.1 provides a general overview of the control scheme in which the controller is employed. Subsection 3.2 formulates the LMPC. Final, Subsection 3.3 describes the learning convergence metric that is used to analyze the convergence of the learning in the controller.

### 3.1 | General overview

As in the previous section, all variables required for the controller formulation are listed in Table 2. Most of these variables serve as input data to the controller and remain fixed during the optimization process. This includes constants such as the gravitational acceleration $(g)$, the train data $(A, B, C, M, L, \tau, a_l, a_f, u_{max}, P_{max})$, the controller settings $(t_s, H_p, H_c, v_{min}, v_{max}, j_{max}, d_{des}, d_{min})$, the line data that depends on the train position $(i, R, v_{DP})$, the information from the preceding train $(s^p$ and $v^p)$, and the information about previous iterations $(v^\lambda, Q_f)$, as explained in the following subsection.

In addition, there are decision variables, denoted by (u),

which are optimized during the solution process. The resulting value of (u) constitutes the control action to be applied.

There are also calculated variables, which correspond to the variables introduced in Section 2 $(\mathbf{X}, j, d, d^{rel})$ and internal variables of the controller $(\varepsilon^{v_{max}}, \varepsilon^{v_{min}}, \varepsilon^{d_{rel}}, \varepsilon^d, \lambda)$. The values of these variables during the optimization step of the control problem depend on the values assigned to the decision variables and must satisfy the controller's constraints.

### 3.2 | LMPC formulation

In this subsection, the LMPC formulation is presented. Since the LMPC relies on solving an optimization problem, we begin by introducing the complete problem formulation in Equation (7). We then provide a detailed explanation of each component, including all constraints and the cost function. This explanation is structured in two parts: the first covers the constraints and cost terms related to the control policy (Equations (8)–(12)), while the second focuses on those associated with the learning process (Equations (13) and (14)).

The LMPC is based on the optimization problem in Equation (7), which depends on the state vector $\mathbf{X}$, the input force decision variable $u$, and the learning variable $\lambda$. This optimization problem is solved for each train $n$ within a convoy of $N$ trains, being $n = 1$ the leader and $n = 2, \dots, N$ the followers, at a time $t$ for the prediction horizon $H_p$, being each solving execution separated by the integration step $t_s$. Therefore, only the first step of the input variable throughout the prediction horizon $(u_{1|t})$ is applied, being recalculated the input variable solutions for future states according to the real evolution of the system at a future time $t$. These executions, which characterize the classical MPCs, ensure a proper control action for the current time $t$, while considering possible future states in the prediction horizon.

Beyond executing the controller at each time step $t$, the controller is also run at every iteration $r$ of the learning process. Each iteration corresponds to a full simulation from the start to the end of the line. As a result, the solution to the optimization problem at a given time step $t$, can differ across iterations due to the incorporation of new data—specifically, data from prior iterations—into the constraint defined by Equation (13). This constraint enables the controller to learn from past experiences, progressively refining the control actions. This iterative learning mechanism distinguishes LMPC from classical MPC. As will be explained later, Equation (13) uses selected data from previous complete simulations as input, encoded through the variable $v^\lambda$. Thus, thanks to the iterations, the LMPC learns from past simulations and improves the behavior of the controlled train.

Hence, the LMPC formulation is based on the optimization problem given in Equation (7).

$$\min_{u_{\cdot|t}} J_T\left(\mathbf{X}_{k|t}, u_{k|t}, \lambda_{k|t}\right) \qquad (7a)$$

subject to:

$$\mathbf{X}_{t|t} = \mathbf{X}(t) \qquad \forall n \in N \qquad (7b)$$

$$\mathbf{X}_{k+1|t} = \mathbf{X}_{k|t} + t_s \dot{\mathbf{X}}_{k|t} \qquad \forall n \in N \qquad (7c)$$

$$u_{k|t} \in \mathcal{U} \qquad \forall n \in N \qquad (7d)$$

$$v_{k|t} \in \mathcal{V} \qquad \forall n \in N \qquad (7e)$$

$$v_{H_p+1|t} \in \mathcal{V}_{H_p+1} \qquad \forall n \in N \qquad (7f)$$

$$d_{k|t}^{rel} \in \mathcal{D} \qquad \forall n \in N\backslash\{1\} \qquad (7g)$$

$$d_{H_p+1|t}^{rel} \in \mathcal{D}_{H_p+1} \qquad \forall n \in N\backslash\{1\} \qquad (7h)$$

$$d_{H_p+1|t} \in \mathcal{D}_{H_p+1} \qquad \forall n \in N\backslash\{1\} \qquad (7i)$$

To begin with, the first part of the optimization problem is explained from this paragraph to Equation (12). This first part involves the constrains and cost terms connected to the control policy.

In Equation (7b), the state vector of the controller at $k = t$ is initialized with the real current state vector as input data. Thus, each controller execution considers the real conditions under which the train is running. Future states are calculated based on the system dynamics for all the prediction horizon, that is, $\forall k \in [t, t + H_p]$, as shown in Equation (7c), which is defined as in Equation (3). These calculated states are used to ensure the fulfillment of constraints in future states and, therefore, to constraint the acceptable solutions of the optimization problem.

In Equation (7d), $\mathcal{U}$ is the set of the allowable values for the input variable $u$. This set consists of the constraints given by Equation (8), in which the control horizon is implemented in Equation (8a), the longitudinal jerk of the Equation (4) is restricted, and the input variable $u_{k|t}$ is constrained by the train's maximum traction and braking forces, in Equation (8c), and power, in Equation (8d). The control horizon specifies the number of states in which the input variable $u$ can have different values throughout the prediction horizon, being constrained to constantly maintain the last calculated value from the state corresponding to the control horizon $H_c$ to the end of the prediction horizon $H_p$. In addition, the longitudinal jerk constraint reduces the oscillations of the decision variable $u$ throughout the control horizon in order to increase comfort. Moreover, Equation (8c) and Equation (8d) constraint the input variable $u$ within the maximum tractive and braking capabilities of the train, according to the tractive and braking effort curves.

$$u_{k+1|t} - u_{k|t} = 0 \qquad \forall k \in [t + H_c, t + H_p - 1] \quad (8a)$$

$$-j_{max} \leq j_{k|t} \leq j_{max} \qquad \forall k \in [t, t + H_p - 1] \qquad (8b)$$

$$-u_{max} \leq u_{k|t} \leq u_{max} \qquad \forall k \in [t, t + H_p] \qquad (8c)$$

$$-P_{max} \leq u_{k|t} \cdot v_{k|t} \leq P_{max} \qquad \forall k \in [t, t + H_p] \qquad (8d)$$

In Equation (7e), $\mathcal{V}$ is the set of the allowable values for the state variable $v$. This set consists of the constraints given by Equation (9), in which the speed is soft constrained between the minimum value $v_{min}$, in Equation (9b), and a speed profile curve $v_{DP}(s_{k|t})$, in Equation (9a), calculated by using a dynamic programming method, as in (Felez, Kim, & Borrelli, 2019). This precomputed speed profile is introduced in the controller as input data and is used as a maximum speed constraint in order to ensure the fulfillment of the line's speed limits according to the traction and braking capabilities of the controlled train. Note also that the speed is constrained to non-negative values in order to avoid reversing as an allowable solution.

$$v_{k|t} \leq v_{DP}(s_{k|t}) + \varepsilon_{k|t}^{v_{max}} \qquad \forall k \in [t, t + H_p] \qquad (9a)$$

$$v_{min} - \varepsilon_{k|t}^{v_{min}} \leq v_{k|t} \qquad \forall k \in [t, t + H_p] \qquad (9b)$$

$$\varepsilon_{k|t}^{v_{min}} \geq 0, \ \varepsilon_{k|t}^{v_{max}} \geq 0 \qquad \forall k \in [t, t + H_p] \qquad (9c)$$

$$v_{k|t} \geq 0 \qquad \forall k \in [t, t + H_p] \qquad (9d)$$

In Equation (7g), $\mathcal{D}$ is the set of the allowable values for the distance between trains $d^{rel}$, as defined in Equation (6). This set consists of the soft constraints given by Equation (10), where the relative braking distance between trains is constrained to remain above the desired distance $d_{des}$. Note also that, unlike previous constraints, this set of constraints is not applied to all the trains within the convoy. As this constraint is the responsible for implementing virtual coupling in the controller and ensuring a safe distance according to this novel signaling system, it is only applied to the followers. In other words, this constraint is not applied to the leader because this train is not virtually coupled to any train running in front of it.

$$d_{des} \leq d_{k|t}^{rel} + \varepsilon_{k|t}^{d_{rel}} \qquad \forall k \in [t, t + H_p] \qquad (10a)$$

$$\varepsilon_{k|t}^{d_{rel}} \geq 0 \qquad \forall k \in [t, t + H_p] \qquad (10b)$$

In Equations (7h) and (7i), $\mathcal{D}_{H_p+1}$ is the terminal set of the allowable values for the distance between trains. This set consists of the soft constraints given by Equation (11), which constraints both the absolute distance between trains $d$ and the relative braking distance between trains $d^{rel}$, as defined in Equations (5) and (6), respectively. In Equation (11a), the constraint involving the relative braking distance between trains is expanded an additional distance, whereas, in Equation (11b), the absolute distance between trains is constrained in order to respect the minimum safe distance $d_{min}$. This Equation (11b) is necessary because there are combinations of speeds which can lead to unsafe negative absolute distances ($d$) if only Equation (11a) is considered. As Equation (7g), note also that Equations (7h) and (7i) are only applied to the followers.

$$d_{des} \leq d_{H_p+1|t}^{rel} + \varepsilon_{H_p+1|t}^{d_{rel}} \tag{11a}$$

$$d_{min} \leq d_{H_p+1|t} + \varepsilon_{H_p+1|t}^{d} \tag{11b}$$

$$\varepsilon_{H_p+1|t}^{d_{rel}} \geq 0, \varepsilon_{H_p+1|t}^{d} \geq 0 \tag{11c}$$

At this point of the description of the optimization problem, the cost function is introduced. The cost function in Equation (7a) consists of two components: the first component is formulated in Equation (12) and involves all the cost terms related to the control policy, whereas the second component, which is formulated in Equation (14), involves the cost term related to the learning process.

As for the first component of the cost function is concerned, Equation (12) considers the cost terms related to the control policy. Apart from considering soft constraints and the jerk, this first component of the cost function defines $v_{DP}(s_{k|t})$ and the $d_{des}$ as a reference for the leader and the follower, respectively. Therefore, the first component (named after $J$) has a different formulation depending on whether the controlled train is a leader or a follower. The difference lies in the considered reference cost term and the inclusion of the distance slack variables in the follower's cost function.

$$J_1 = \sum_{k=1}^{k=H_p-1} \left(\frac{j_{k|t}}{j_{max}}\right)^2 + \sum_{k=1}^{k=H_p} \left(\frac{\varepsilon_{k|t}^{v_{min}}}{v_{max}} + \frac{\varepsilon_{k|t}^{v_{max}}}{v_{max}}\right) \tag{12a}$$

$$J_d = \sum_{k=1}^{k=H_p} \left(\frac{d_{k+1|t}-d_{des}}{d_{des}}\right)^2 + \sum_{k=1}^{k=H_p+1} \left(\frac{\varepsilon_{k|t}^{d_{rel}}}{d_{des}}\right) + \frac{\varepsilon_{H_p+1|t}^{d}}{d_{des}} \tag{12b}$$

$$J_v = \sum_{k=1}^{k=H_p-1} \left(\frac{v_{k+1|t}-v_{DP}(s_{k+1|t})}{v_{max}}\right)^2 \tag{12c}$$

$$J(\mathbf{X}_{k|t}, \mathbf{u}_{k|t}) = \begin{cases} J_1 + J_v & n = 1 \\ J_1 + J_d & n > 1 \end{cases} \tag{12d}$$

Finally, the second part of the optimization problem is explained from this paragraph to Equation (14). This second part involves the constrains and cost term connected to the learning process.

In Equation (7f), $\mathcal{V}_{H_p+1}$ denotes the terminal set defining the allowable values of the state variable $v$ at the end of the prediction horizon. This terminal set is constructed based on the constraint defined in Equation (13), which embeds the learning mechanism into the controller while also guaranteeing system stability (Rosolia & Borrelli, 2019). The core idea behind the learning mechanism is to guide the system toward an intermediate speed—computed relative to the actual speeds achieved in previous iterations over a horizon-length interval and given the train's current position. This terminal set is updated at each iteration by incorporating data from past iterations, introduced to the LMPC through the variable $v^\lambda$.

$$v_{H_p+1|t} = v_{\cdot|t}^{\lambda}(s_{k|t}) * \lambda_{\cdot|t} \qquad \forall r \in [1, r] \tag{13a}$$

$$\lambda_{\cdot|t} \geq 0 \qquad\qquad\qquad \forall r \in [1, r] \tag{13b}$$

$$\sum \lambda_{\cdot|t} = 1 \qquad\qquad\qquad \forall r \in [1, r] \tag{13c}$$

As introduced at the beginning of the problem formulation, the variable $v^\lambda$ receives a subset of the data of each complete previous simulation as input data to the controller. This subset consists of the running speed for the current position of the train in previous iterations, plus the running speeds for the following known states in a horizon time window with a $H_p + 1$ length in each iteration. Therefore, $v^\lambda$ is a column vector with 1 row and $r \cdot (H_p + 1)$ columns and $\lambda$ is a row vector of decision variables with $r \cdot (H_p + 1)$ rows and 1 column, which considers the costs associated to the real speeds in previous iterations by means of the learning cost $Q_f$ in the column vector of 1 row and $r \cdot (H_p + 1)$ columns. This learning cost $Q_f$ has been defined with identical terms as the specified for $J(\mathbf{X}_{k|t}, \mathbf{u}_{k|t})$ in Equation (12) and evaluated with the simulation results of the previous iterations. Hence, the second component of the cost function ($J_L$) is defined in Equation (14).

$$J_L(\lambda_{k|t}) = Q_f(s_{k|t}) * \lambda \tag{14}$$

Therefore, the final and complete cost function of the optimization problem that is minimized in Equation (7a) consists of the two components formulated in Equation (12) and the learning component introduced in Equation (14). This final and complete cost function ($J_T$) is defined in Equation (15).

$$J_T(\mathbf{X}_{k|t}, \mathbf{u}_{k|t}, \lambda_{k|t}) = J(\mathbf{X}_{k|t}, \mathbf{u}_{k|t}) + J_L(\lambda_{k|t}) \tag{15}$$

### 3.3 | Learning convergence metric

In this subsection, a learning convergence metric is defined in order to analyze the convergence of the LMPC. The metric, which is provided in Equation (16), is based on the controller cost given by Equation (12), but is evaluated for the actual states and forces obtained during simulation.

$$L_1 = \left(\frac{j(t)}{j_{max}}\right)^2 + \frac{\varepsilon^{v_{min}}(t)}{v_{max}} + \frac{\varepsilon^{v_{max}}(t)}{v_{max}} \tag{16a}$$

$$L_d = \left(\frac{d(t)-d_{des}}{d_{des}}\right)^2 + \frac{\varepsilon^{d_{rel}}(t)}{d_{des}} + \frac{\varepsilon^{d}(t)}{d_{des}} \tag{16b}$$

$$L_v = \left(\frac{v(t)-v_{DP}(s(t))}{v_{max}}\right)^2 \tag{16c}$$

$$L(\mathbf{X}_{k|t}, \mathbf{u}_{k|t}) = \begin{cases} L_1 + L_v & n = 1 \\ L_1 + L_d & n > 1 \end{cases} \tag{16d}$$

Thus, when the total average cost throughout the simulation at each iteration, calculated as the average of the individual costs at each instant, which is given by Equation (16), converges to a steady value, the LMPC is considered to have converged. This steady value can be studied by means of the reduction in the variation of this total average cost in one iteration compared to the previous iteration, as shown in the figures of the simulation results.

# 4 | SIMULATIONS AND RESULTS

In this section, we conduct simulations in order to assess the controller presented in this paper.

The simulations are performed in three different line categories for passenger trains: a metro line, a regional train, and a high-speed train. The profile lines for each category are summarized in Figure 1, while the parameters of the corresponding trains used for each category are listed in Appendix A. These line profiles and train parameters have been extracted from real reports published and publicly made available by the Spanish railway authorities and the railway companies, respectively.

In each line category, the convoy considered for simulation consists of two trains: the leader (train $n = 1$) and the follower (train $n = 2$). This configuration allows the study of a real implementation in which only two trains can simultaneously stop at each station.

Note also that each line can be divided in inter-stop segments. Thus, the complete metro line simulations comprise three segments; the regional train, two; and the high-speed line, one.

For each one of these line categories, two simulations are performed. The first simulation tests the behavior of the LMPC when applied only to the follower, with the leader using a classical MPC, as in (Vaquero-Serrano & Felez, 2023a). The second simulation tests the behavior of the LMPC when simultaneously applied to the leader and the follower. Therefore, in the first simulation, the follower is the only train that performs the learning procedure, while, in the second simulation, both the leader and the follower perform the training. The results are compared with a classical MPC as the one presented in (Vaquero-Serrano & Felez, 2023a) in both simulations.

All simulations have been developed with Yalmip (Lofberg, 2004), MATLAB, and a computer with an i7-1365U 1.8 GHz processor and 32 GB of RAM. Due to the fact that loading data from previous iterations is necessary, the simulation loop presents the form shown in the pseudocode contained in Algorithm 1.

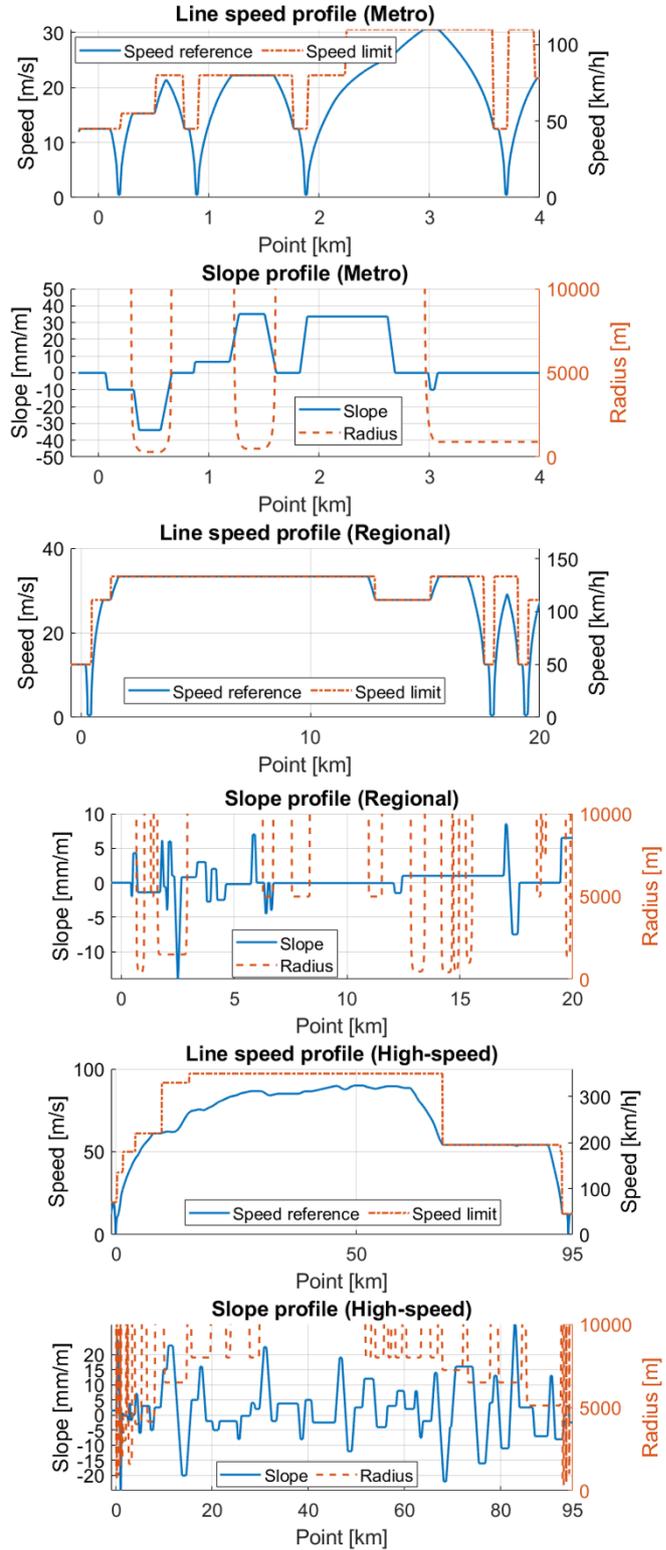

**FIGURE 1** Profile line: speed references, speed limits, slopes and radii of the different simulated lines (metro, regional and high-speed).

## 4.1 | Simulation 1: LMPC applied only to the follower

Simulation 1 assesses the behavior of the LMPC when applied only to the follower, with the leader using a classical MPC. The results are firstly analyzed in terms of the LMPC behavior on its own, and then, compared with a classical MPC as the one presented in (Vaquero-Serrano &

Felez, 2023a). The simulation has been done in the three-line categories, although the figures only show one of them for simplicity, whereas Table 3 shows the complete results for all the train categories and for the complete given line and its line segments.

### ALGORITHM 1 LMPC implementation

1:  Initialize algorithm: load train and line data, and initialize the controller;
2:  **for** r = 0 to 10 **do**
3:      Initialize simulation ($t = 0$);
4:      Load data from previous simulations (iterations) and calculate $Q_f$ as in Equation (12) for each final simulated state in those simulations;
5:      **for** $t = 0$ to line's end **do**
6:          Retrieve data according to current position: retrieve line data and data from previous simulations (iterations);
7:          Retrieve available data of the preceding train;
8:          **for** $n = 1$ to $N$ **do**
9:              **if** LMPC applicable to train $n$ **do**
10:                  Execute LMPC: Solve optimization problem in Equation (7) and obtain $u$;
11:             **else**
12:                  Execute classical controller and obtain $u$;
13:             **end if**
14:         **end for**
15:         Apply $u_{1|t}$ to the simulated dynamics in Equation (1) of each train and obtain the train states for the next simulation step;
16:     **end for**
17: **end for**

First, the LMPC on its own preserves good VC conditions, as seen in Figure 2, in which the regional results are shown as an example. The VC conditions consist of the results in terms of distance between trains, speeds, and accelerations and decelerations. These VC conditions are considered to be good (or acceptable) if the accelerations and decelerations remain within the traction and braking limits and if the integrity of the convoy is maintained by means of speeds and a close distance between trains. Regarding the speeds, the convoy is considered to maintain its integrity when all the trains that are part of the convoy can arrive at the station and stop approximately at the same time. This result can be seen at stops in the time-speed plots. Regarding the distance between trains, outside the stations, the distance between trains must increase as the speed increases in order to maintain safe conditions according to virtual coupling constraints.

All the above characteristics of the VC conditions can be observed in Figure 2, in which the regional results are shown as an example. Note that, when analyzing the LMPC on its own, the results of the figure that are being analyzed correspond uniquely to the 'Follower (Iteration 10)' plot line. In the time-speed plot (Figure 2a), it is shown that all the trains that are part of the convoy arrive at the station and stop approximately at the same time. Outside the stations, as the speed increases, the distance between trains (Figure 2c) also increases to maintain safe conditions according to virtual coupling constraints. Conversely, when the speed reduces, the distance between trains also reduces according to the virtual coupling constraints. Furthermore, the accelerations and decelerations remain within the traction and braking limits, as shown in Figure 2b. In Figure 2, note also that the leader is also represented at the beginning and the end of the iterations associated with the follower in order to show that the leader's behavior is not modified during the learning process.

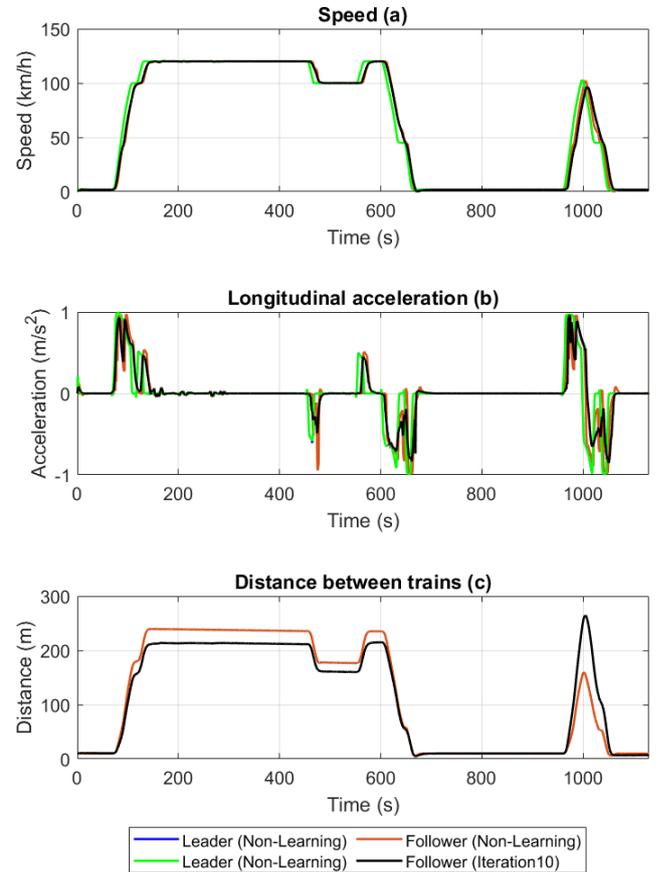

**FIGURE 2** Results of the regional train for the LMPC applied only to the follower. (a) speed, (b) longitudinal acceleration, and (c) distance between trains.

Second, compared with a classical MPC, the LMPC reduces energy consumption of the follower, without significantly affecting the virtual coupling conditions, as shown in Table 3. In this table, note that the maximum reached distance between trains of the complete line considers the difference, in terms of maximum reached distance, between the classical MPC and the LMPC throughout the complete line and, therefore, it allows the comparison of different line segments when the maximum reached distance changes the line segment in which it presents. In general, as seen in the table, the variation of the maximum reached distance of the complete line equals the value of one of the line segments (for instance, in the metro line, the variation of +0.9% corresponds to the segment 3, in which the maximum distance between trains in the line is reached), but there is a noticeable exception in the regional train, as assessed in Section 4.3. This exception is due to the fact that the maximum reach distance in the complete regional line changes from the first segment to the second one, as respectively seen in the 'Follower (Non-Learning)' and 'Follower (Iteration 10)' plot lines in Figure 2c, where 'Non-Learning' denotes the classical MPC and 'Iteration 10' denotes the LMPC finally implemented once the controller has learnt. Therefore, note that, when comparing with the classical MPC, the differences between the 'Non-Learning' and the 'Iteration 10' plot lines are being analyzed.

The energy consumption reduction is achieved thanks to a better optimization of the control forces than in MPC, which manifests in an average acceleration/force reduction. In fact, this optimization of the control forces is particularly important in negative gradients, where the intrinsic acceleration of the profile line is used to reduce traction forces to maintain the same speed, as can be seen in the first segment of the metro line.

Regarding also the energy consumption reduction, note that there is an important difference between the different train categories, which can be summarized in two conclusions, being the first a consequence of the second. First, the results reveal that the specific energy consumption reduction in the regional train is lower than in the metro train but it is greater than in the high-speed train. Second, the obtained results show that the largest specific energy consumption reductions are achieved in line segments with multiple, frequent, and significant variations of the speed, by means of accelerations and decelerations, not generally maintaining a constant speed for long running distances. This occurs in all the segments of the metro line and in the second segment of the regional line. On the contrary, the lowest specific energy consumption reductions occur in line segments in which the speed remains approximately constant for long running distances. This occurs in the high-speed line and in the first segment of the regional line. Hereinafter, this situation will be referred to as a constant-speed segment.

In order to exemplify the previous paragraph, the different segments of the regional line are analyzed. It can be seen that the energy consumption reduction in the first constant-speed segment is lower than in the second, which is characterized by an acceleration and a deceleration without a long constant-speed segment. The first segment is characteristic of long-distance lines, where the distance between stops is long enough to maintain a constant speed for long periods of time. Conversely, the second segment is similar to what happens on metro lines, whose main characteristic is that the distance between stops is not usually long enough to maintain a constant speed for long periods of time, which causes multiple accelerations and decelerations. In other words, the first segment is similar to a high-speed line in terms of constant-speed segments, whereas the second segment is closer to a metro line, in which speed variations occur more frequently.

As a matter of the coherence of the results, note also that the sum of the variations of the energy consumption of the different segments equals the variations of the complete line.

Finally, experience shows that LMPC converges between the fifth and tenth iteration, as shown in Figure 3 for the Metro as an example.

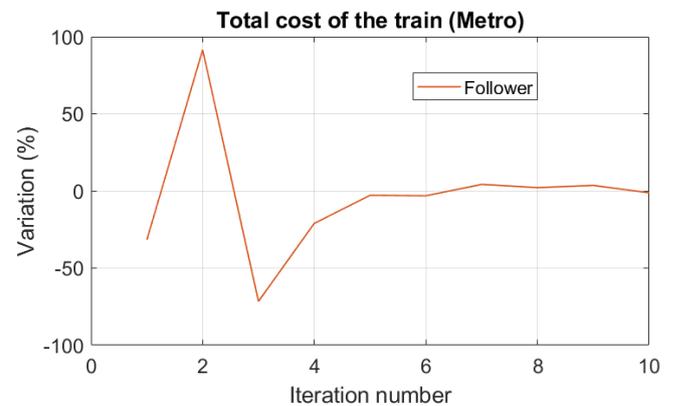

**FIGURE 3** Learning convergence for the Metro train when the LMPC is applied only to the follower.

## 4.2 | Simulation 2: LMPC applied to the leader and the follower

Simulation 2 assesses the behavior of the LMPC when applied to both the leader and the follower. As in the previous simulation, the results are firstly analyzed in terms of the LMPC behavior on its own ('Iteration 10' plot line) and, secondly, they are compared with a classical MPC ('Non-Learning' plot line). Simultaneously, the results should be analyzed differentiating between leader and follower.

Regarding the LMPC applied to the leader, the LMPC-controlled leader on its own (train $n = 1$) achieves the desired speeds according to the DP profile, as shown in Figure 4a, in which the regional results are shown as an example. In addition, the results reveal an acceptable behavior in terms of accelerations and decelerations, which remain within the traction and braking limits (Figure 4b).

Compared with a classical MPC, the LMPC reduces energy consumption of the leader, as shown in Table 4. This energy consumption reduction is achieved thanks to a better optimization of the control forces than in MPC, which manifests in an average traction force reduction. As in Simulation 1, this optimization of the control forces is particularly important in negative gradients. Therefore, the proposed LMPC manages to reduce the energy consumption of each train, leading to an energy consumption reduction of the entire convoy. In addition, due to the fact that LMPC reduces energy consumption in each train, the energy consumption reduction is greater in Simulation 2 than in Simulation 1, in which the LMPC is not applied to the leader.

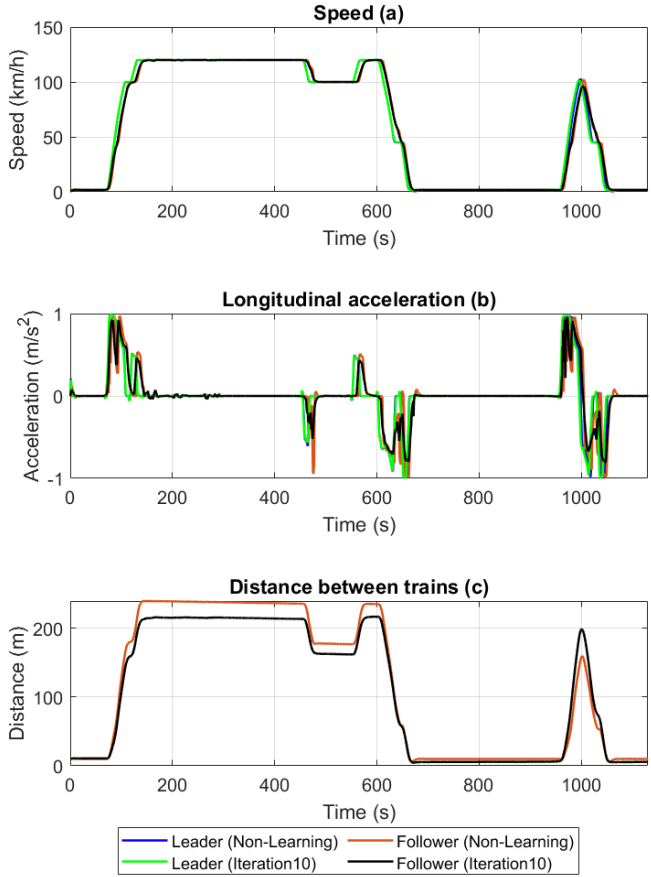

**FIGURE 4** Results of the regional train for the LMPC applied to the leader and the follower. (a) speed, (b) longitudinal acceleration, and (c) distance between trains.

Likewise, regarding the energy consumption reduction, note that there is also an important difference between the different train categories. As in Simulation 1, the results reveal that it can be concluded that the specific energy consumption reduction of the proposed LMPC, in both the leader and the follower, is larger for segments that resemble a metro line and lower for segments that resemble a high-speed line, in terms of the distances between stops and the number of accelerations and decelerations involved.

Regarding the LMPC applied to the follower (train $n = 2$), the results lead to the same conclusions as in Simulation 1, which means that the application of the LMPC to the leader does not negatively affect the general behavior of the LMPC in the follower (Figure 4). In fact, note that the results for the followers are similar in both Simulation 1 and Simulation 2.

Despite indicating the same general behavior and tendency, note that the specific numerical results are not the same as in Simulation 1 because the running conditions of the in-front train have varied due to the application of the LMPC to the leader. Therefore, the LMPC for the follower reveals promising results for real applications in which the running conditions of the leader can vary due to uncertainties and disturbances.

In Table 4, as in Table 3, the maximum reached distance between trains on the complete line considers the difference, in terms of maximum reached distance, between the classical MPC and the LMPC throughout the complete line. In this table, it can be seen that there is no variation of the line segment in which the maximum reached distance happens for the regional train, albeit the variations present the same tendency as in Simulation 1, as seen in Figure 4c. As in Simulation 1, the reason for this behavior is also assessed in Section 4.3.

As a matter of coherence of the results, as in Simulation 1, note also that the sum of the variations of the energy consumption of the different segments is equal to the variations of the complete line. In addition, the sum of the variations of the different trains equals the variations of the convoy of both the different segments and the complete line.

Finally, as in Simulation 1, experience shows that LMPC converges between the fifth and tenth iteration in both: the leader and the follower, as shown in Figure 5, in which the high-speed train is shown as an example.

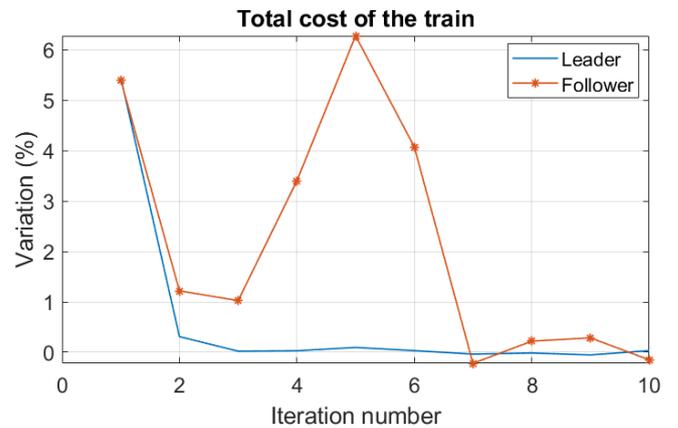

**FIGURE 5** Learning convergence of the high-speed train.

TABLE 3 Simulation 1 results: variation of variables with respect to a classical MPC when the LMPC is applied only to the follower.

| Parameter variation | Segment line | Metro | Regional | High-speed |
|---|---|---|---|---|
| Maximum reached distance between trains (m) | 1 | 5.1 (6.7%) | -24.6 (-10.3%) | 161 (10%) |
| | 2 | 3.1 (4.4%) | 105.1 (65.9%) | |
| | 3 | 0.9 (0.9%) | | |
| | **Complete** | **0.9 (0.9%)** | **25.3 (10.6%)** | **161 (10%)** |
| Force average (kN) | Complete | -1.8 (-17.8%) | -0.4 (-3.6%) | -1.2 (-2.4%) |
| Specific energy (kJ/t*km) | 1 | -47.1 (-22.5%) | -2.6 (-3.3%) | -1.7 (-1.1%) |
| | 2 | -19.0 (-4.9%) | -38.2 (-12.5%) | |
| | 3 | -13.0 (-3.1%) | | |
| | **Complete** | **-79.1 (-7.8%)** | **-40.8 (-10.6%)** | **-1.7 (-1.1%)** |

TABLE 4 Simulation 2 results: variation of variables with respect to a classical MPC when the LMPC is applied to the leader and follower.

| Parameter variation | Train | Segment line | Metro | Regional | High-speed |
|---|---|---|---|---|---|
| Maximum reached distance between trains (m) | Follower | 1 | 0.5 (0.6%) | -22.7 (-9.5%) | 160 (9.9%) |
| | Follower | 2 | 3.0 (4.3%) | 39.1 (24.5%) | |
| | Follower | 3 | 0.1 (0.1%) | | |
| | **Follower** | **Complete** | **0.1 (0.1%)** | **-22.7 (-9.5%)** | **160 (9.9%)** |
| Force average (kN) | Leader | Complete | -1.7 (-11.2%) | 0.0 (-0.4%) | 0.7 (1.5%) |
| | Follower | Complete | -1.6 (-15.8%) | -0.3 (-2.8%) | -1.1 (-2.2%) |
| Specific energy (kJ/t*km) | Leader | 1 | -25.8 (-11.5%) | -0.0 (0.0%) | -0.4 (-0.2%) |
| | Leader | 2 | -6.3 (-1.6%) | -6.3 (-2.1%) | |
| | Leader | 3 | -2.2 (-0.5%) | | |
| | **Leader** | **Complete** | **-34.3 (-3.3%)** | **-6.3 (-1.6%)** | **-0.4 (-0.2%)** |
| | Follower | 1 | -46.9 (-22.4%) | -2.6 (-3.3%) | -1.8 (-1.1%) |
| | Follower | 2 | -19.8 (-5.1%) | -40.1 (-13.1%) | |
| | Follower | 3 | -13.1 (-3.1%) | | |
| | **Follower** | **Complete** | **-79.8 (-7.8%)** | **-42.7 (-11.1%)** | **-1.8 (-1.1%)** |
| | Convoy | 1 | -72.7 (-16.8%) | -2.6 (-1.7%) | -2.2 (-0.7%) |
| | Convoy | 2 | -26.1 (-3.3%) | -46.5 (-7.6%) | |
| | Convoy | 3 | -15.3 (-1.8%) | | |
| | **Convoy** | **Complete** | **-114.1 (-5.5%)** | **-49.1 (-6.4%)** | **-2.2 (-0.7%)** |

## 4.3 | Comparison of Simulation 1 and 2 for further follower results. The multi-objective problem

Some observations about the results of the regional train can be made with respect to other trains. In both Simulations 1 and 2, it can be observed that the variation of the maximum reached distance between the regional trains reduces when the constant speed segment is analyzed (segment R1), whereas it increases in the second segment (segment R2), which is closer to a metro line in terms of varying-speed segments.

In segment R1, apart from reducing the energy consumption by 3.3% for the follower, the LMPC manages to reduce the distance between trains by 9.5% in Simulation 2. However, this does not mean that the LMPC will reduce the distance between trains in all constant speed segments to which it might be applied. Instead, this result outlines that, apart from always reducing the energy consumption, the LMPC will try to improve other terms of the cost function depending on each particular situation, train and line characteristics. This improvement will be a compromise among the different competitive, and sometimes opposing, objectives of the cost function. For instance, in this segment R1, reducing the distance between trains is possible because the distance reference cost can be reduced without increasing the costs associated with the distance and speed constraints.

In contrast, there is another example of these competitive objectives that manifests the opposing effect: the increase of the maximum distance between trains in the high-speed train. As seen in Figure 1, the simulated high-speed line is formed by two main subsegments within the unique constant speed segment between stops: a 350 km/h segment (350-speed subsegment), which involves a constant speed that the simulated train cannot maintain due to high slopes and a relatively low maximum power of the train, and a 200 km/h segment (200-speed subsegment), which involves a speed that the simulated train can effectively maintain. Due to the fact that virtual coupling imposes longer distances for higher speeds, the maximum distance variation shown in Tables 3 and 4 for the high-speed train corresponds to the 350-speed subsegment, in which this variation is positive for Simulations 1 and 2, meaning that the distance between trains increases. Similarly, and in contrast to the regional train, the distance between trains also increases in the 200-speed segment. This situation is due to the fact that the cost term related to the distance constraint is optimized further to the detriment of the distance reference cost. As a result, in this subsegment, the maximum distance between trains increases 97 m (8.5%) in Simulation 1 and 98 m (8.6%) in Simulation 2.

In segment R2, there is another particularly interesting example of these competitive objectives, whose effects can significantly increase the distance between trains, as seen in Simulation 1. This increase is due to the fact that the LMPC tries to improve the cost terms related to the distance constraints in a segment with a high uphill followed by an important downhill. In this situation, the energy consumption and the constraint satisfaction lead to a large increment in the distance between trains in Simulation 1. However, this increment is lower in Simulation 2 because the leader also improves its energy consumption in this particular segment, reducing the speed differences between both trains and, therefore, achieving a distance increase between the classical MPC and Simulation 1.

## 4.4 | Time of computation

After the introduction of the additional constraints and costs terms through Equations (13) and (14), respectively, there could be concerns about the computational time needed to solve the optimization problem.

In Table 5, these concerns are solved by analyzing the time of computation needed in Simulations 1 and 2 (Simul). In this table, the total time of computation of the LMPC ($t_{LMPC}$) for each train category throughout the simulation is compared with the time of computation of the classical MPC ($t_{MPC}$), revealing that the time of computation of the LMPC is larger than in the classical MPC. In addition, there is not a clear criterion to differentiate the time of computation obtained for each train category. However, when compared with the time that is being simulated in each case ($t_{simul}$), results show that the time of computation is lower than the time that is being simulated and remains within the same order of magnitude.

Therefore, despite being more complex than a classical MPC, the LMPC can be still implemented in real time.

**TABLE 5** Time of computation.

| Category | Simul | Train | $t_{LMPC}$ (s) | $t_{MPC}$ (s) | $t_{simul}$ (s) |
|---|---|---|---|---|---|
| Metro | 1 | Follower | 119.4 | 81.3 | 341.0 |
| | 2 | Leader | 146.4 | 68.0 | 341.0 |
| | | Follower | 109.0 | 81.3 | 341.0 |
| Regional | 1 | Follower | 454.7 | 372.3 | 1430.0 |
| | 2 | Leader | 432.6 | 351.3 | 1430.0 |
| | | Follower | 402.4 | 372.3 | 1430.0 |
| High-speed | 1 | Follower | 498.4 | 449.8 | 1847.2 |
| | 2 | Leader | 550.5 | 539.5 | 1847.2 |
| | | Follower | 533.1 | 449.8 | 1847.2 |

## 5 | CONCLUSIONS

In this paper, a LMPC that improves the control policy through a terminal set and that is applied in the Autonomous Driving and Train Control rail domain for virtual coupling has been proposed. The objective of this controller has been the optimization of the behavior of the virtually coupled train convoy, while maintaining the control policy of an equivalent classical MPC without learning capacities. This control policy involves running the trains at the maximum speed allowed by the limitations and characteristics of the line in accordance with the characteristics of the train that is being controlled. In addition, the control policy involves keeping the different trains in the convoy as close as possible, while maintaining safe conditions and avoiding the risk of a collision. Thus, the LMPC has been firstly trained by repeatedly running the train throughout the line in what have been denominated as 'iterations.' Once the LMPC has undergone several iterations, it contains within its formulation the experience from these previous iterations and can be directly applied to control a virtually coupled train and improve its behavior.

The proposed LMPC has been simulated for three categories of trains: a metro, a regional train, and a high-speed train. The simulations considered the possibility of applying the LMPC only to the follower, as well as of applying the LMPC to both: the leader and the follower.

The main results conclude that the LMPC achieved a significant reduction in energy consumption in the trains in which it is applied, without affecting the virtual coupling and maintaining similar operating conditions, with similar travel times and speeds. This energy consumption reduction was achieved thanks to a better optimization of the control forces in LMPC than in MPC and it was larger the shorter the distance between stops and the more acceleration and deceleration segments involved. Moreover, experience showed that the proposed LMPC converges between the fifth and tenth iterations. Furthermore, the time of computation of the LMPC remains sufficiently close to the time of computation of a classical MPC for a real implementation.

Future research will focus on the integration of the LMPC with other advances in virtual coupling, such as (Luo, Tang, Chai, & Liu, 2024), for an industrial application. This industrial application will allow testing the LMPC in a real environment and the analysis of the eventual effect of the introduction of real data to the LMPC instead of the data obtained from the application of the LMPC in a simulation loop.

## 6 | REFERENCES

# 7 | APPENDIX A: TRAIN PARAMETERS

**TABLE A.1** Train parameters.

| Parameter | Metro | Regional | High-speed | Parameter | Metro | Regional | High-speed |
|---|---|---|---|---|---|---|---|
| M (kg) | 99.972e+3 | 247.48e+3 | 457.4e+3 | $t_s$ (s) | 0.2 | 0.2 | 0.2 |
| L (m) | 54.9 | 107.36 | 200 | $H_p$ | 20 | 20 | 20 |
| A (N) | 1216.13 | 1804.5 | 3383.5 | $H_c$ | 2 | 1 | 1 |
| B (N/(m/s)) | 117.39 | 68.87 | 114.55 | $v_{max}$ (m/s) | 30.6 | 69.4 | 97.2 |
| C (N/(m/s)²) | 2.97 | 4.91 | 7.32 | $v_{min}$ (m/s) | 0.5 | 0.5 | 0.5 |
| $\tau_d, \tau_b, \tau$ (s) | 0.7 | 0.7 | 0.7 | $jerk_{max}$ (m/s³) | 0.98 | 0.98 | 0.98 |
| $F_1^{max}$ (N) | 97972.56 | 242.53e+3 | 143.44e+3 | $d_{des}$ (m) | 10 | 10 | 10 |
| $F_2^{max}$ (N) | 150e+3 | 242.53e+3 | 224.126e+3 | $d_{min}$ (m) | 4/6 | 4/6 | 4/6 |
| $P_1^{max}$ (W) | 1.584e+6 | 4e+6 | 8.8e+6 | $a_l$ (m/s²) | 1.25 | 1.25 | 0.625 |
| $P_2^{max}$ (W) | 1.584e+6 | 4e+6 | 8.8e+6 | $a_f$ (m/s²) | 1 | 1 | 0.5 |
| $\mu$ | 0.15 | 0.1 | 0.05 | | | | |